\documentclass[12pt]{article}
\usepackage{latexsym}

\newcommand{\vspc}{\vspace{.15 in}}

\usepackage {latexsym}

\begin{document}
   \title{ Shifted Vertex Operator Algebras}
\author{Chongying Dong\footnote{Supported by NSF grants, a China NSF 
grant and  a Faculty
research grant from  the University of
California at
Santa Cruz} and Geoffrey Mason\footnote{Research supported by the 
NSF and UC Santa Cruz} \\
Department of
Mathematics \\
University of California, Santa Cruz.}
\date{}
   \maketitle

\vspc
\noindent
{\sc Abstract:} We study the properties of \emph{shifted
vertex operator algebras}, which are vertex algebras derived from a
given theory by
shifting the conformal vector. In this way, we are able to exhibit
large numbers of vertex operator algebras which are regular
(rational and $C_2$-cofinite) and yet are pathological in one
way or another.

\vspc

\vspc
\section{Introduction}
\setcounter{equation}{0}
The motivation to develop the results in the present paper
stems from a desire to have available a source of
examples of vertex operator algebras which can be used to
test ideas. The most familiar types of vertex operator algebras
tend to have a number of features in common. A large portion
   of the
literature (in both mathematics and physics)
  concerned with vertex
operator algebras $V$ assumes that
$V$ is of \emph{CFT-type}, that is the underlying Fock space has the shape
\begin{equation}\label{eq: CFT-type}
V = \mathbf{C1} \oplus V_2 \oplus V_3 ...
\end{equation}
Another common assumption is that $V$ is \emph{self-dual}
  in the
sense that the dual module $V'$ is isomorphic to the
  adjoint module
$V$.  Deviation from any of these assumptions could be considered
somewhat pathological.

  \vspc
   It is  a natural question to ask
whether there are significant classes of vertex operator algebras
exhibiting such pathologies.
   We are particularly interested in the
case in which the vertex operator algebra $V$
   is  nevertheless well-behaved from the representation-theoretic
point-of-view. By this
   we generally have in mind that $V$ is
\emph{rational} in the sense of \cite{DLM3},
   that is,  the module
category $V-Mod$ is semisimple.
   Closely related to this is the condition that $V$ is \emph{regular}
in the sense of \cite{DLM1}. By
   \cite{ABD}, this is the same as requiring that $V$
is both rational and $C_2$-cofinite (cf. \cite{Z}). A case of 
particular interest
involves \emph{holomorphic} vertex operator algebras,
   which by definition are rational and for which the adjoint module
is the unique simple module.

\vspc
    Consider, then, the following pathologies that
a vertex operator algebra $V$ might exhibit:
   \begin{eqnarray*}
&&(i) \ V \ \mbox{has nonzero negative weight spaces;} \\
&&(ii) \ V \ \mbox{is not self-dual as $V$-module;} \\
&&(iii)\ \mbox{the zero weight space is degenerate: dim} V_0 \geq 2.
   \end{eqnarray*}
It is easily seen (see below for the easy proof) that (i)
$\Rightarrow$ (iii),  so that there are
six remaining possible combinations of (i)-(iii) which may occur. We will show
in Section 5 that
there are regular, simple vertex operator algebras of all six types.
Indeed, it will be clear from our
construction that such vertex operator algebras \emph{exist in profusion}.

\vspc
   A related issue that arises is concerned with the question of the
\emph{finiteness}
of the number of vertex operator algebras satisfying various
conditions. Among other things, we construct examples of the following:
\begin{eqnarray*}
&&(1) \mbox{ infinitely many nonisomorphic regular, simple vertex operator} \\
   && \mbox{ algebras  with the \emph{same} partition function} \
Z_V(q) = \mbox{Tr}_V q^{L(0) - c/24} \\
&&\mbox{ and \emph{equivalent} module categories;} \\
&&(2) \mbox{ for any integer $c$ divisible by $8$, infinitely many
nonisomorphic } \\
&&\mbox{  holomorphic vertex operator algebras of central charge $c$}.
\end{eqnarray*}
One may arrange for the vertex operator algebras in (1) to be
holomorphic, for example,
or that none of them are self-dual.
Thus there is no finiteness theorem for holomorphic vertex operator
algebras, say, of a given central charge, that is analagous to
Minkowski's finiteness
theorem for positive-definite self-dual lattices.
   Our examples \emph{are} consistent with the possibility of a
finiteness theorem for
(rational, say) vertex operator algebras with given values of both
central charge $c$ and \emph{effective} central charge $\tilde{c}$
(cf. \cite{DM1}).

\vspc
The method of construction for all of these examples is not difficult. It
involves taking a vertex operator algebra with known properties (we
use lattice theories (\cite{LL}) as they are particularly accessible)
and \emph{shifting} the Virasoro element while at the same time retaining
   the overall space of vertex operators. For this reason, we
call this general class of theories
\emph{shifted vertex operator algebras}. This idea was already used
in \cite{DLinM} to a limited extent, and by Matsuo and Nagatomo \cite{MN},
who considered shifting the conformal vector in
the Heisenberg vertex operator algebra.
   As far as we know however, there is presently no systematic
discussion of the properties
of shifted vertex operator algebras in the literature That is what we
carry out here.

   \vspc
   As we have suggested, once the basic
construction is presented, the implementation
   in explicit examples
is not difficult. Nevertheless, the type of pathological vertex
operators
   that can result seem not to be well-known.
   We discuss the case of lattice theories set in Section 2, showing
how a shifted lattice
theory gives rise to a particular type of vertex algebra that we call
a \emph{$\mathbf{C}$-graded vertex operator algebra},
or more generally a $k$-graded vertex operator algebra, $k$  a
subgroup of $\mathbf{C}$. Vertex operator algebras in the usual sense
correspond to the case $k = \mathbf{Z}$.
Sections 3-5 are devoted to elucidating the properties of
shifted lattice theories in the case $k = \mathbf{Z}$, and among 
other things we
construct the various examples alluded to above. In the final Section
6, we show how
our constructions can be put in a more abstract context, so that one
can define shifted vertex operator algebras in a quite general setting.

\vspc
\section{$\mathbf{C}$-graded Vertex Operator Algebras}
\setcounter{equation}{0}

We begin with the

\vspc
\noindent
{\sc Definition:} A \textit{\mbox{$\mathbf{C}$}-graded vertex
operator algebra} is a
quadruple $(V, Y, \mathbf{1}, \omega)$ consisting of a $\mathbf{C}$-graded
complex linear space (the \emph{Fock space})
\begin{equation}\label{eq: linspac}
V = \bigoplus_{ r \in \mathbf{C}}V_r,
\end{equation}
a linear map
\begin{eqnarray}\label{eq: Ymap}
Y: &V& \rightarrow  ({\rm End} V)[[z,z^{-1}]], \nonumber \\
    &v& \mapsto Y(v, z) = \sum_{n \in \mathbf{Z}}v_n z^{-n-1}
\end{eqnarray}
together with a pair of distinguished states $\mathbf{1}, \omega \in V$
(the \emph{vacuum} and \emph{Virasoro} elements respectively). The
following axioms are imposed:\\

(a) For any $u,v\in V,$ $u_nv=0$ if $n$ is sufficiently. \\

(b) Grading: Each homogeneous subspace $V_r$ has finite dimension, and
the $\mathbf{C}$-grading is truncated from below in the following sense:
the implication
\begin{eqnarray}\label{eq: grade}
V_r \neq 0 \Rightarrow \mbox{Re}(r) \geq | \mbox{Im}(r)|
   \end {eqnarray}
holds for \textit{all but finitely many} $r$. (Here, Re($r$) and
Im($r$) refer to real and imaginary
parts of the complex number $r$.)\\

(c) Locality: For any pair of states $u, v \in V,$ there exits a nonnegative
integer $n$ such that
$$(z_1-z_2)^n[Y(u, z_1),Y(v,z_2)]=0.$$

(d) Creativity: If $v \in V$ then
\begin{eqnarray*}
Y(v, z)\mathbf{1} = v + \sum_{n < -1}v(n)\mathbf{1}z^{-n-1}.
\end{eqnarray*}

(e) Conformality: The state $\omega$ is a Virasoro element. That is,
if $Y(\omega, z) = \sum L(n)z^{-n-2}$
then there is a complex number $c$ (central charge) such that
\begin{equation}\label{eq: Virasororeln}
[L(m), L(n)] = (m-n)L(m+n) + \frac{(m^3-m)}{12}\delta_{m, -n}c Id.
\end{equation}
Moreover
\begin{eqnarray*}
\frac{d}{dz}Y(v, z) = Y(L(-1)v, z)
\end{eqnarray*}
for $v \in V$, and $L(0)$ is a semisimple operator such that
\begin{eqnarray*}
V_r = \{ v \in V | L(0)v = rv \}.
\end{eqnarray*}

A $\mathbf{C}$-graded vertex operator algebra is therefore a vertex
algebra with a particular kind
of grading induced by the $L(0)$ operator. Note that from (\ref {eq:
grade}), there are only finitely many
values of $r$ for which $V_r \neq 0$ and Re($r) < 0$. In particular, if
the grading is \emph{real} (i.e. $V_r \neq 0 \Rightarrow r \in
\mathbf{R}$), then the
truncation condition simply says that $V_r = 0$ for all small enough
$r$. If the grading is \emph{integral}
(i.e. $V_r \neq 0 \Rightarrow r \in \mathbf{Z}$), then $V$ is just a vertex
operator algebra in
the usual sense. For an additive subgroup $k \subseteq \mathbf{C}$,
we say that $V$ is
\emph{$k$-graded} in case $V$ is $\mathbf{C}$-graded and $V_r \neq 0
\Rightarrow r \in k$.

\section{Shifted Lattice Theories}
\setcounter{equation}{0}

In this Section we construct a large number of $\mathbf{C}$-graded
vertex operator algebras by the method of \emph{shifting}
the Virasoro element of a lattice vertex operator algebra $V_L$.
   We will see that this idea leads to a rich
source of examples of $\mathbf{C}$-graded vertex operator algebras. The basic
set-up, which we discuss below in somee detail, is presented in
Section 4 of \cite{DLinM}.

\vspc
   Let $L$ be a positive-definite even lattice of rank $l$, with inner product
\begin{eqnarray*}
(\ , \ ): L \ \mbox{x} \ L \rightarrow \mathbf{Z}.
\end{eqnarray*}
Let $H  = \mathbf{C} \otimes L$
be the corresponding complex linear space equipped with the $\mathbf{C}$-linear
extension of $( \ , \ )$. In fact, it is useful to introduce
some related spaces, and for this purpose we recall that
   the $\mathbf{Z}$-dual, or \emph{dual lattice} of $L$, is defined by
\begin{eqnarray*}
L^o = \{ f \in \mathbf{R} \otimes L \  | \  (f, \alpha) \in \mathbf{Z}, \
\mbox{all} \ \alpha \in L \}.
\end{eqnarray*}
Then for a (nonzero) additive subgroup $k \subseteq
\mathbf{C}$, set $H_k = k \otimes L^o \subseteq H$.

\vspc
Let $(M, Y, \mathbf{1},
\omega_L)$ be the free bosonic vertex operator based on $H$, and \\
$(V_L, Y, \mathbf{1}, \omega_L)$ the
corresponding lattice vertex operator algebra. Both of these theories
have central charge $l$. The
Fock space for $V_L$ is
\begin{equation}\label{eq: latticefockspace}
V_L = M(1) \otimes \mathbf{C}[L]
\end{equation}
where $\mathbf{C}[L]$ is the group algebra of $L$. For more
background we   refer the reader to \cite{FLM}.

   \vspc
For a state $h \in H \subseteq V_L$ we set
\begin{equation}\label{eq: shiftedvirasoro}
\omega_h = \omega_L + h(-2)\mathbf{1}.
\end{equation}
Note that $\omega_h \in (V_L)_2$. We are going to consider the quadruple
\begin{equation}\label{eq: shiftedvir}
(V_L, Y, \mathbf{1}, \omega_h),
\end{equation}
which we denote by $V_{L, h}$ when it is convenient.
This means that we are looking at the Fock space (\ref{eq: latticefockspace})
equipped with the \emph{same} set of fields as $V_L$ and the \emph{same}
   vacuum state. However
the Virasoro state has been shifted, and consequently the conformal structure
and grading are modified.

\vspc
   {\sc Theorem 3.1}: The quadruple (\ref {eq: shiftedvir}) is a
$\mathbf{C}$-graded
vertex operator algebra with central charge
\begin{equation}\label{eq: shiftedcc}
c_h = l - 12(h, h).
\end{equation}
If $h \in H_k$ for an additive subgroup $k \subseteq \mathbf{C}$,
then $V_{L, h}$ is $k$-graded. In particular, if $h \in L^o$
then $V_{L, h}$ is a vertex operator algebra.

\vspc
   {\sc Proof:} Note that parts of this result are already contained in
\cite{DLinM} and \cite{MN}.
   Since $V_L$ is a vertex operator algebra, and in view of
the fact that the Fock space, fields and vacuum state of $V_{L, h}$
are identical to those for $V_L$, it is evident
that the locality and creativity axioms for $V_{L, h}$ hold. So we only need
check the Virasoro and grading axioms.

\vspc
   That $\omega_h$ is a Virasoro element is well-known (loc. cit.) We run
through the details to get the central charge. There are the
following relations:
\begin{equation}\label{eq: heisenbergrelns}
[h(m), h(n)] = m \delta_{m, -n}(h, h) \mbox{Id},
\end{equation}
\begin{equation}\label{eq: primaryfld}
L(m)h = 0 \ \mbox{for} \ m \geq 1, \ L(0)h = h.
\end{equation}
Moreover,
\begin{eqnarray*}
(h(-2)\mathbf{1})(n) = (L(-1)h)(n) = -nh(n-1).
\end{eqnarray*}
So if we set
\begin{equation}\label{eq: shiftedomegaops}
Y(\omega_h, z) = \sum_{n \in \mathbf{Z}} L_h(n)z^{-n-2},
\end{equation}
then
\begin{equation}\label{eq: shiftedLs}
L_h(n) = L(n) - (n+1)h(n).
\end{equation}
Now using (\ref {eq: Virasororeln}), (\ref{eq: heisenbergrelns}) and
(\ref{eq: primaryfld})
we calculate that
\begin{eqnarray*}
[L_h(m), L_h(n)] = (m-n)L_h(m+n) + \frac{(m^3-m)}{12}\delta_{m, -n}(l
- 12(h, h))\mbox{Id}.
\end{eqnarray*}
So the central charge is indeed $l - 12(h, h)$. It is an important
feature of (\ref{eq: shiftedLs})
that
\begin{eqnarray*}
L_h(-1) = L(-1),
\end{eqnarray*}
from which it is clear that the derivative property also holds in $V_{L, h}$.

\vspc
   We turn to consideration of the grading. Obviously from (\ref {eq:
shiftedLs}) we have
\begin{eqnarray*}
L_h(0) = L(0) - h(0).
\end{eqnarray*}
Now $h(0)$ is a semisimple operator on $V_L$ with action
\begin{eqnarray*}
h(0) : u \otimes e^{\alpha} \mapsto (h, \alpha) u \otimes e^{\alpha}
\end{eqnarray*}
for $u \in M(1), \alpha \in L$. Then clearly $L_h(0)$ is also semisimple
and if $u  \in M(1)_n$ then
\begin{eqnarray*}
L_h(0) : u \otimes e^{\alpha} \mapsto (n - (h, \alpha)) u \otimes e^{\alpha}.
\end{eqnarray*}
So the eigenspace $(V_{L, h})_r$ for $L_h(0)$ with eigenvalue $r$ is spanned by
states $u\otimes e^{\alpha}$ with $u \in M(1)_n$ and $n \geq 0$, $\alpha
\in L$, and satisfying
\begin{equation}\label{eq: eigenstates}
    n + \frac{1}{2}(\alpha, \alpha) - (h, \alpha) = r .
\end{equation}
Let us write
\begin{equation}\label{eq: realineq}
r = x + iy, \; x, y \in \mathbf{R},
\end{equation}
\begin{equation}\label{eq: cmpxineq}
h = a + ib, \;  a, b \in \mathbf{R}\otimes L.
\end{equation}
Then (\ref {eq: eigenstates}) tells us that
\begin{equation}\label{eq: realandim}
n + \frac{1}{2}(\alpha, \alpha) - (a, \alpha) = x, \\
-(b, \alpha) = y.
\end{equation}
Now because $( \ , )$ is positive-definite, all but finitely many
$\alpha \in L$ satisfy $(\alpha, \alpha) \geq 4(a \pm b, a\pm b)$.
In this case, an application of the Schwarz inequality leads to
\begin{eqnarray*}
&&(a, \alpha) \pm (b, \alpha) = (a \pm b, \alpha)
               \leq |(a \pm b, \alpha)| \\
               &\leq& \sqrt{(a \pm b, a \pm b)(\alpha, \alpha)}
               \leq \frac{1}{2}(\alpha, \alpha),
\end{eqnarray*}
so that
\begin{eqnarray*}
|y| = |(b, \alpha)|  \leq \frac{1}{2}(\alpha, \alpha) - (a, \alpha) =
x-n \leq x.
\end{eqnarray*}

This proves that $V_{L, h}$ satisfies the truncation condition (\ref
{eq: grade}).
In order to complete the proof of Theorem 3.1, it remains to show that
each of the eigenspaces $(V_{L, h})_r$ has finite dimension.
But this follows from what we have already done. Indeed
the first equation of (\ref {eq: realandim}) can be written in the form
\begin{equation}\label{eq: newvareq}
\frac{1}{2}(\alpha - a, \alpha - a) = x - n + \frac{1}{2}(a, a).
\end{equation}
Because $( , )$ is positive-definite and $L$ discrete,
for fixed $x, a$ there
are only finitely many choices for $\alpha \in L$ and $0 \leq n \in \mathbf{Z}$
which satisfy (\ref {eq: newvareq}). Finally, if
$h \in H_k$ for an additive subgroup $k \subseteq \mathbf{C}$,
   it follows from (\ref {eq: eigenstates}) that $V_{L, h}$
is $k$-graded. This completes the proof of the Theorem. $\Box$

\vspc
   We refer to the $\mathbf{C}$-graded vertex operator algebras
$V_{L, h}$ as \emph{shifted} vertex operator algebras.

\vspc
   There is a more precise approach to the grading on
the vertex operator algebras $V_{L, h}$, for which we consider
the corresponding \emph{partition functions}.
For a $\mathbf{C}-graded$ vertex operator algebra (\ref {eq: linspac}) of
central charge $c$, we define the partition
function to be the usual formal $q$-expansion
\begin{equation}\label{eq: partfunc}
Z_V(q) = \mbox{Tr}_V q^{L(0) - c/24}.
\end{equation}
For example, it is well-known that we have
\begin{equation}\label{eq: lattpartfunc}
Z_{V_L}(q) = \frac{\theta_L(q)}{\eta(q)^l},
\end{equation}
where $\eta(q)$ is Dedekind's \emph{eta function}
\begin{eqnarray*}
\eta(q) = q^{1/24}\prod_{n=1}^\infty(1 - q^n),
\end{eqnarray*}
and $\theta_L(q)$ is the \emph{theta function} of $L$.
For any coset $C = L - h \subseteq H$ of $L$ in $H$ we define Fock spaces
\begin{eqnarray*}
\mathbf{C}[L-h] = \oplus_{\alpha \in L} \mathbf{C}e^{\alpha - h},
\end{eqnarray*}
\begin{eqnarray*}
V_{L-h} = M(1) \otimes \mathbf{C}[L-h],
\end{eqnarray*}
   and the formal sum
\begin{eqnarray*}
\theta_C(q) &=& \sum_{f \in C}q^{(f, f)/2} = \sum_{\alpha \in L}
q^{(\alpha-h, \alpha-h)/2} \\
&=& q^{(h, h)/2}\sum_{\alpha \in L}q^{(\alpha, \alpha)/2 - (h, \alpha)}.
\end{eqnarray*}
$\theta_L(q)$ is just the case when $C = L$. Of particular interest
are the cosets $L - \lambda$ which are contained in the \emph{dual
lattice} $L^o$.
   We know from \cite{D1} that for $\lambda \in L^o$, the Fock space
   $V_{L- \lambda}$  naturally carries the structure of a simple $V_L$-module.
Moreover as $L - \lambda$ ranges over the elements of $L^o/L$,
we obtain in this way each simple $V_L$-module exactly once.
One defines the partition function
of  the Fock spaces $V_{L-h}$-modules by the obvious extension of
(\ref {eq: partfunc}),
and we have (generalizing (\ref {eq: lattpartfunc}))
\begin{equation}\label{eq: shiftpartfunc}
Z_{V_{L-h}}(q) = \frac{\theta_{L-h}(q)}{\eta(q)^l}.
\end{equation}

Let us now consider the partition function for the $\mathbf{C}$-graded
vertex operator algebra $V_{L, h}$. We see from Theorem 3.1 and its proof that
\begin{eqnarray*}
Z_{V_{L,h}}(q) &=&  \mbox{Tr}_{V_{L,h}} q^{L_h(0) - c_h/24}\\
&=& \mbox{Tr}_{V_{L}} q^{L(0) - h(0) - l/24 + (h, h)/2}\\
&=& Z_M(q) \mbox{Tr}_{\mathbf{C}[L]}q^{L(0) - h(0) + (h, h)/2}\\
&=& Z_M(q) \sum_{\alpha \in L} q^{(\alpha-h, \alpha-h)/2}\\
&=&\frac{\theta_{L-h}(q)}{\eta(q)^l}.
   \end{eqnarray*}
Using (\ref {eq: shiftpartfunc}), we have proved

\vspc
   {\sc Proposition 3.2:} The partition function of the shifted lattice vertex
operator algebra satisfies
\begin{equation}\label{shiftpartfuncidentity}
Z_{V_{L,h}}(q) = Z_{V_{L-h}}(q).
\end{equation} $\Box$

\vspc
This result has several interesting corollaries. Suppose
first that $h \in L^o$.
Then $V_{L,h}$ is a ($\mathbf{Z}$-graded)
vertex operator algebra by Theorem 3.1. On the other hand, by Proposition
3.2 and the discussion that preceded it, the partition function of $V_{L, h}$
is identical to that of the simple $V_L$-module $V_{L-h}$. It is now
easy to establish

\vspc
{\sc Theorem 3.3:} Let $L$ be a (nonzero) positive-definite even
lattice of rank $l$,
and let $N$ be any simple $V_L$-module. Then there
are infinitely many pairwise nonisomorphic ($\mathbf{Z}$-graded)
vertex operator algebras,
\emph{all} of which have the \emph{same} partition function as $N$.

\vspc
{\sc Proof:} Let $N = V_{L-\lambda}$ for some $\lambda \in L^o$.
We have only to consider the vertex operator algebras $V_{L,
h}$ where $h$ ranges over an infinite sequence of elements in the coset
$L - \lambda$
such that $(h, h)$ is strictly increasing. By (\ref {eq: shiftedcc}),
no two of these vertex operator algebras
are isomorphic, since the corresponding central charges are strictly
decreasing. On the
other hand, we have just explained that each of them has
partition function equal to $Z_N(q)$. $\Box$

\section{Modules over $V_{L, h}$}
\setcounter{equation}{0}

   In this Section we continue the discussion of the modules and
partition functions
for the vertex operators $V_{L, h} \ (h \in L^o$) begun in $\S3$.
Recall from $\S3$
that for a vertex operator algebra
$V$ and $V$-module $N$ the partition function is
\begin{eqnarray*}
Z_{N,V}(q) = Z_N(q) = Tr_N q^{L(0) - c/24}.
\end{eqnarray*}
For what follows, we refer the reader to \cite{DLM2} for results and 
terminology
concerning $V$-modules.

\vspc
   {\sc Theorem 4.1:} Let $h \in L^o$. There are equivalences between
the categories
   of weak, admissible and ordinary $V_L$-modules, and the categories
of weak, admissible and ordinary $V_{L, h}$-modules, respectively. In
particular, $V_{L, h}$ is a \emph{regular}
vertex operator algebra.

\vspc
   {\sc Proof:} Suppose that $(N, Y_N)$ is any weak , admissible or
ordinary module for $V_L$. Thus
for $v \in V_L, Y_N(v, z)$ is the field on $N$ determined by $v$. We turn
$N$ into a (weak) $V_{L, h}$-module in the obvious way. Namely, take the
\emph{same} fields $Y_N(v, z)$, now with $v$ considered as a state in
$V_{L. h}$.
It is clear from the definition of module (over $V_L$) that in this way,
$N$ indeed becomes a weak $V_{L, h}$-module. In the same way, because $V_L$ and
$V_{L, h}$ share the same set of fields, a weak $V_{L, h}$-module
is also a $V_L$-module. In this way, the equivalence of the category
of weak modules
for $V_L$ and $V_{L, h}$ follows.

\vspc
   Now take $N$ to be a simple $V_L$-module, so that $N = V_{L-\lambda}$
for some $\lambda \in L^o$. Then certainly $N$ is an irreducible module
for the operators spanned by Fourier modes of states in $V_{L, h}$. Moreover,
the argument used to establish Proposition 3.2 shows that
\begin{eqnarray*}
Z_{N, V_{L,  h}}(q) = Z_{V_{L-\lambda-h}}(q) =
\frac{\theta_{L-\lambda-h}}{\eta(q)^l}.
\end{eqnarray*}
In particular, $N$ is indeed a module over $V_{L, h}$. Finally, let
$N$ be a simple
module for $V_{L, h}$. Then $N$ is a simple weak module
for $V_L$, and  by Theorem 3.16 of \cite{DLM2}, $N$ is therefore a 
simple (ordinary)
module for $V_L$. All parts of the Theorem now follow easily from what we have
established, together with the result \cite{DLM2} that $V_L$ is
a regular vertex operator algebra. $\Box$

\vspc
   From the proof of Theorem 4.1, we see that the object map of the
categorical equivalence
\begin{equation}\label{eq: functor}
F : V_L-Mod \longrightarrow V_{L, h}-Mod
\end{equation}
   is just an identification $F(N) = N$.

\vspc
   Suppose that
$L = L^o$ is self-dual. This is equivalent \cite{D1} to the assertion 
that $V_L$ has
a \emph{unique} simple module, namely $V_L$ itself. Indeed $V_L$
is \emph{holomorphic}. That is, it is regular and has a unique simple module.
By Theorem 4.1, we can conclude that $V_{L, h}$ is also a holomorphic
vertex operator algebra whenever $h \in L$.

\vspc
   {\sc Theorem 4.2:}. Let $c$ be any integer divisible $8$. Then there
are infinitely
many pairwise nonisomorphic, holomorphic, $C_2$-cofinite vertex
operator algebras of central charge $c$.

\vspc
\noindent
{\sc Remark:} By results from \cite{Z} and \cite{DLM3}, we know that if $V$ is
holomorphic
and $C_2$-cofinite, then we necessarily have $8|c$.

\vspc
   {\sc Proof:}  Let $d$ be a positive integer satisfying
\begin{eqnarray*}
24d + c > 0.
\end{eqnarray*}
For each integer $r \geq d$, let $(L_r, h_r)$ be a pair consisting of
a self-dual, even lattice $L_r$ of rank $24r + c$ and an element $h _r \in L_r$
which satisfies $(h_r, h_r) = 2r$. Such a sequence of pairs is easily
constructed:
for example, we can take $L_r$ to be the orthogonal direct sum of
$3r+\frac{c}{8}$
copies of the $E_8$ root lattice. From the discussion preceding the
statement of
Theorem 4.2, it follows that each of the shifted vertex operator algebras
$V_{L_r, h_r}$ is holomorphic, and from (\ref{eq: shiftedcc}) we see that
the central charge of each $V_{L_r, h_r}$ is $c$.

\vspc
   Next we assert that if $r \neq s$ then $V_{L_r, h_r}$ and $V_{L_s, h_s}$
are \emph{not} isomorphic as vertex operator algebras. Indeed, by
(\ref {shiftpartfuncidentity}),
the partition function of $V_{L_r, h_r}$ coincides with that of
$V_{L_r}$ itself,
and hence has the shape
\begin{eqnarray*}
Z_{V_{L_r, h_r}}(q) = q^{-r-c/24}(1 + \cdots).
\end{eqnarray*}
Our assertion follows from this. Finally, the $C_2$-cofiniteness condition
for a vertex operator algebra $V$ refers to the finiteness of the
codimension of the
subspace $C_2(V)$ of $V$ spanned by states $u(-2)v$ for $u, v \in V$.
In view of the fact
that $V_L$ and $V_{L, h}$ share the same states and fields,
it is evident that we have $C_2(V_L) = C_2(V_{L, h})$. Since one knows (
\cite{DLM3}, Proposition 12.5) that $V_L$ is $C_2$-cofinite, then so too are
the vertex operator algebras $V_{L, h}$. This completes the proof of
the Theorem.
$\Box$

\vspc
   Next we consider the question of \emph{duality} of modules. Recall 
\cite{B}, \cite{FHL}
that if $V$ is a vertex operator algebra and $(N, Y_N)$ a $V$-module,
then the restricted
dual $N'$ of $N$ becomes a $V$-module if we inflict upon it the
fields $Y_N'(v, z)$ \emph{adjoint}
to $Y_N(v, z)$. Precisely, the adjoint is defined for $v \in V$ using
\begin{eqnarray*}
Y_N(v, z)^{\dagger} = Y_N(e^{zL(1)}(-z^{-2})^{L(0)}, z^{-1}).
\end{eqnarray*}
Then $(N', Y_N')$ is a $V$-module, called the \emph{dual} module of $N$.
For $f \in N', v \in V, n \in N$ and $< \ , >: N' \otimes N
\rightarrow \mathbf{C}$ the canonical pairing, we have
\begin{eqnarray*}
<Y_N'(v, z)f, n> = <f, Y_N(v, z)^{\dagger}n>.
\end{eqnarray*}
We say that $N$ is \emph{self-dual} in case there is an
isomorphism of $V$-modules $N \cong N'$. This applies, in particular,
to the adjoint module $V$ itself. For example,
a holomorphic vertex operator algebra is necessarily self-dual.

\vspc
   {\sc Proposition 4.3:} The equivalence of categories (\ref{eq: functor})
preserves dualities. That is, we have $F(N') = F(N)'$ for any $V_L$-module
$N$.

\vspc
   {\sc Proof:} We have already seen that any $V_L$-module is \emph{ipso facto}
a $V_{L, h}$-module with the same set of fields. It follows from what we have
said that the adjoint operators are also the same, and the
Proposition follows immediately.
$\Box$

\vspc
   Now we know that for a lattice vertex operator algebra $V_L$,
the dual of the simple $V_L$-module $V_{L - \lambda}$ ($\lambda \in L^o$)
is $V_{L + \lambda}$. Hence, the same is true if we regard this
pair as modules over a shifted vertex operator algebra $V_{L, h}$, $h \in L^o$.
   Then by Propositions 3.2 and 4.3, the partition
functions for $V_{L, h}$ and its dual satisfy
\begin{eqnarray*}
Z_{V_{L, h}}(q) = Z_{V_{L-h}}(q), \\
Z_{V_{L, h}'}(q) = Z_{V_{L+h}}(q).
\end{eqnarray*}
  From this and Theorem 4.1,  it follows that the following is true:

\vspc
   {\sc Theorem 4.4:} Let $L$ be an even lattice and $h \in L^o$. Then
the shifted vertex operator algebra $V_{L, h}$ is self-dual if, and only if,
$2h \in L$. $\Box$

\section{Types of simple vertex operator algebras}
\setcounter{equation}{0}

In this Section we consider various types of simple vertex operator
algebras. The
attributes of $V$ that concern us here are related to the nature of the maps
in the following commutative diagram:
   \begin{eqnarray*}
             \begin{array}{ccccc}
&&V_{-1}&\stackrel{L(1)}{\longleftarrow}&V_0 \\
&&\downarrow&&\downarrow \\
   V_{-1}& \stackrel{L(1)}{\longleftarrow} &V_0
&\stackrel{L(1)}{\longleftarrow}& V_1 \\
   \downarrow& & \downarrow& & \\
V_0&\stackrel{L(1)}{\longleftarrow}&V_1&&
    \end{array}
\end{eqnarray*}
   where all vertical maps are $L(-1)$.

\vspc
   {\sc Lema 5.1}\cite{L}: Let $V$ be a simple vertex operator algebra. Then
\begin{equation}\label{eq: codim}
\mbox{dim} \ V_0/L(1)V_1 \leq 1.
\end{equation}

   {\sc Proof:} Li proved \cite{L} that for
a vertex operator algebra $V$, the space Hom$_V(V, V')$
of $V$-module maps of the adjoint module $V$ into the dual
module $V'$ has dimension equal to that of $V_0/L(1)V_1$. But
if $V$ is also simple then by Schur's Lemma, the Hom space
in question has dimension at most $1$. The Lemma is proved. $\Box$

\vspc
   {\sc Lemma 5.2}: Suppose that $V$ is a vertex operator algebra. Then
\begin{eqnarray*}
   \mbox{dim} V_0 > \mbox{dim} \ V_{-1}.
\end{eqnarray*}
   In particular, if $V_0 = \mathbf{C}\mathbf{1}$ then $V$ has no
nonzero negative weight spaces.

\vspc
{\sc Proof:}  We know \cite{DLinM} that $L(-1): V_n \rightarrow V_{n+1}$
is an injection
   as long as $n \neq 0$. It therefore suffices to
show that $\mathbf{1} \notin L(-1)V_{-1}$. Suppose
   that $v \in
V_{-1}$ satisfies $L(-1)v = \mathbf{1}$. Bearing in mind that $L(1)$
and $L(-1)$
   generate an algebra $S \cong sl_2$ of operators on V,
we see that $v$ generates an
   $S$-submodule for which $\mathbf{1}$
is the \emph{unique} highest weight vector
   (up to scalars). The
structure theory of $sl_2$-modules shows that this is not possible.

$\Box$

\vspc
\noindent
{\sc Definition:} Suppose that $V$ is a simple vertex operator
algebra. We say that $V$ has \emph{type
$I$} if $V_0 = L(1)V_1$ and \emph{type II} if $L(1)V_1 \neq V_0$;
\emph{type A}
if dim$V_0 > 1$ and \emph{type B} if dim$V_0 = 1$; \emph{type +} if
$V_{-1} = \{0\}$
and \emph{type -} if $V_{-1} \neq \{0\}$. In other words, $V$ has
type $II$ or type $I$
according to whether it is self-dual or not; it has type $B$ or type
$A$ according to whether
the vacuum is nondegenerate or not; and type $-$ or $+$ according to whether
it has nonzero negative weight spaces or not.

\vspc
One can combine these qualities to obtain eight types
of simple vertex operator algebras in all.
Not all of them exist, however.  It is a consequence of Lemma 5.2 that
there can be \emph{no} simple vertex operator algebras of types $IB-$
or $IIB-$.
We present a table of the possibilities.  Lemma 5.1 shows that the
only two possibilities for
the first column are $0$ and $1$.

\vspc
\noindent
\( \begin{array}{cccccc}
\hspace{2.5 cm} \underline{ \mbox{codim}L(1)V_1}&
\underline{\mbox{dim}V_{-1}}& \underline{\mbox{dim}V_0}&
\underline{\mbox{Type}}&
   \underline{\mbox{Exist?}}
\\
   \hspace{3 cm}0&>0&>1&IA-& \mbox{Yes} \\
\hspace{3 cm}0& >0&1&IB-&\mbox{No}\\
\hspace{3 cm}0&0&>1&IA+& \mbox{Yes} \\
\hspace{3 cm}0&0&1&IB+& \mbox{Yes}\\
\hspace{3 cm}1&>0&>1&IIA-& \mbox{Yes}\\
\hspace{3 cm}1&>0&1&IIB-& \mbox{No}\\
\hspace{3 cm}1&0&>1&IIA+& \mbox{Yes}\\
\hspace{3 cm}1&0&1&IIB+&\mbox{Yes}
\end{array} \)

\vspc
\noindent
The "existence" column indicates whether a given type exists or not.
We have already explained
nonexistence for types $IB-$ and $IIB-$.
   We next
   give some further explicit constructions of shifted lattice
vertex operator algebras which will establish the existence of examples of
each of the remaining six types. From what we have already proved in
previous Sections,
all of our examples are simple, $C_2$-cofinite, regular vertex
operator algebras.

\vspc
\noindent
{\sc Example 1:} Let $0\leq k \leq 2N$ be integers with $N \geq 2$.
Let $L = \mathbf{Z}\alpha$
be the $1$-dimensional lattice spanned by $\alpha$ where $(\alpha,
\alpha) = 2N$, and
take $h = k\alpha/2N \in L^o$. Thus $V_{L, h}$ is a vertex operator
algebra. We will discuss
the properties of $V_{L, h}$ as $k$ varies within the indicated range.

\vspc
\noindent
(i) $k = 0$. Type $IIB-$. Here $V_{L, h}$ is simple the lattice
theory $V_L$, which is
self-dual, has no nonzero negative weight spaces, and has nondegenerate vacuum.

\vspc
\noindent
(ii) $0<k<N$. Type $IB+$. In this case $2h \notin L$, so $V_{L, h}$
is not self-dual
by Theorem 4.4. If $\beta = m\alpha \in L$ and $u \in H$, then when
considered as
a state in $V_{L, h}$, we have
\begin{eqnarray*}
wt(u \otimes e^{\beta}) = wt(u) + m(mN-k).
\end{eqnarray*}
This follows from (\ref {eq: eigenstates}). Since $0<k<N$ then
$m(mN-k) \leq 0$ if, and only if, $m=0$. It follows that the vacuum space of
$V_{L, h}$ is nondegenerate, so that $V_{L, h}$ has the type stated.

\vspc
There are
three other essentially different cases, each of which can be analyzed
in the same way. We merely state the result in each case:

\vspc
\noindent
(iii) $k=N$. Type $IIA+$.

\vspc
\noindent
(iv) $N<k<2N$. Type $IA-$.

   \vspc
\noindent
(v) $k=2N$. Type $IIA-$.

\vspc
In this way we have constructed 5 of the six possible types of
vertex operator algebras. Only type $IA+$ remains
to be accounted for. We get this type in the next example.

\vspc
\noindent
{\sc Example 2:} We take $L$ to be the root lattice of type $A_l$
with $l \geq 2$. Let $h \in L^o$ be such that $L + h$ generates $L^o/L$.
As this latter quotient has order $l+1 \geq 3$ then $V_{L, h}$
is not self-dual by Theorem 4.4. A calculation only slightly more
complicated that the one above shows that $V_{L, h}$ has
no nonzero negative weight spaces, that is it has type $+$.
Moreover the zero weight space is spanned by states of the form
$\mathbf{1} \otimes e^{\beta}$ where $\beta \in L$ is either
$\mathbf{1}$ or a connected sum
\begin{eqnarray*}
\beta = \alpha_1 + \alpha_{l-1} + ... + \alpha_t, 1\leq t \leq l,
\end{eqnarray*}
where $\alpha_1, ...\alpha_l$ is a fundamental system of roots.
   So $V_{L, h}$ has the asserted type $IA+$, indeed we have shown that the
zero weight space has dimension $l+1$.
It follows that for all integers $n\geq 1$, there are (simple) vertex operator
algebras $V$ of type $IA+$ for which the zero weight space $V_0$
has dimension $n$

\section{An Abstract Approach}
\setcounter{equation}{0}

We have made use of Proposition 3.2 mainly in the case that $h \in L^o$.
In general it says that the partition function of the shifted vertex operator
algebra $V_{L, h}$ coincides with that of a twisted $V_L$-module,
namely the one determined by the coset $L - h$ in $H/L$. (See \cite{DM2}
for twisted modules over lattice theories in the case of automorphisms
of finite order.) In this Section we show how to
extend this to a general setting. The main idea is to utilize a
certain operator $\Delta( \ , z)$ introduced by Li \cite{L}, and which we
have used elsewhere \cite{DLM4}.

\vspc
We work with the following set-up: $V$ is a vertex operator algebra
of central charge $c$ and
   $h \in V_1$ satisfies the following conditions:
\begin{eqnarray*}
(i)&&h \ \mbox{is a primary state, i.e.} \ L(n)h = 0, n \geq  1; \\
(ii)&&h(0) \ \mbox{is semisimple with real eigenvalues}; \\
(iii) &&h(n)h = 0 \ \mbox{for} \  0 \leq n \neq 1,  \mbox{and} \ 
h(1)h \in \mathbf{C1}; \\
(iv)&&[h(m), h(n)] = m \delta_{m, -n}h(1)h \mbox{Id}.
\end{eqnarray*}
The same argument
given during the proof of Theorem 3.1 shows that
$\omega_h = \omega - h(-2) \mathbf{1}$ is
a Virasoro element in $V$ of central charge
\begin{eqnarray*}
c_h = c - 12h(1)h,
\end{eqnarray*}
where we identify $h(1)h$ with the scalar multiple of $\mathbf{1}$
to which it is equal. We continue to use
the notation (\ref {eq: shiftedomegaops}), whence (\ref {eq:
shiftedLs}) still holds.

\vspc
\noindent
Now set
\begin{eqnarray*}
\Delta(h, z) = z^{h(0)}\mbox{exp}\{ - \sum_{k \geq 1} \frac{h(k)}{k}
(-z)^{-k}\}.
\end{eqnarray*}

\vspc
{\sc Proposition 6.1:} Suppose that $(M, Y_M)$ is a
$V$-module. Then
   $(M, Y_ {M, h})$ is a $e^{- 2 \pi i h(0)}$-twisted
$V$-module, where we set
\begin{eqnarray*}
Y_{M, h}(v, z) = Y_M(\Delta(h, z)v, z) \ \mbox{for} \ v \in V.   \ \ \  \Box
\end{eqnarray*}

\vspc
If the eigenvalues of $h(0)$ are \emph{rational},
  this result has been
proved in \cite{L}. But the same proof works in our more general
setting if we use the definition of twisted module given in \cite{DLinM}
for an automorphism whose eigenvalues lie on the unit circle.
\vspc

We will prove an analog of Proposition 3.2, namely
\vspc

\noindent
{\sc Proposition 6.2:} We have
\begin{eqnarray*}
Z_{(M, Y_ {M, -h})}(q) = \mbox{Tr}_Mq^{L_h(0) - c_h/24}.
\end{eqnarray*}

\vspc
{\sc Proof:} The reader may verify that
\begin{eqnarray*}
\Delta( - h, z)\omega = \omega - hz^{-1} + \frac{1}{2}h(1)hz^{-2},
\end{eqnarray*}
so that the corresponding zero mode operator is
\begin{eqnarray*}
L_{\Delta, - h}(0) = L(0) - h(0) + \frac{1}{2}h(1)h.
\end{eqnarray*}
It follows that
\begin{eqnarray*}
Z_{(M, Y_ {M, -h})}(q) &=& \mbox{Tr}_Mq^{L_{\Delta, -h}(0) - c/24} \\
                                            &=& \mbox{Tr}_M q^{ L(0) -
h(0) + \frac{1}{2}h(1)h - c/24} \\
                                            &=& \mbox{Tr}_M q^{L_h(0) - c_h/24},
\end{eqnarray*}
as required.    $\Box$

\vspc
Now assume that $h(0)$ has \emph{integral} eigenvalues, so that $e^{2
\pi i h(0)}$
is the \emph{trivial} automorphism of $V$. Then $(M, Y_{M, h})$ is
a $V$-module by Proposition 6.1, and in particular it is $\mathbf{Z}$-graded.
Taking $V = M$, Proposition 6.2 together with the argument used in the
proof of Theorem 3.1 show that $(V, Y, \omega_h, \mathbf{1})$
is a vertex operator algebra of central charge $c_h$. We denote it by
$V_h$. Evidently, $V_h$ is
a type of shifted vertex operator algebra that  includes
the shifted lattice theories considered earlier. Furthermore,
we see as in the proof of Theorem 4.1 that $V$ and $V_h$
have equivalent module categories (of various types). For example,
if $V$ is rational then so too is $V_h$, and the set of partition
functions of the simple modules for the two
vertex operator algebras agree. In particular, in this case the effective
central charge is an invariant, i.e. it is the same for each vertex
operator algebra.

We end the paper with a natural question. Having constructed vertex 
operator algebras exhibiting various pathologies by modifying the 
conformal vector of a `nice' vertex operator algebra, one can ask 
whether all such pathologies arise in this way.

\vspc
\noindent
{\sc Question:} Suppose $V=(V,Y,{\bf 1},\omega)$ is a simple 
$\mathbf{Z}$-graded  rational vertex operator algebra. Is it true 
that there exists $h\in V_1$ satisfying conditions (i)-(iv) 
 above, 
and such that $(V,Y,{\bf 1}, \omega_h)$ is of CFT type?

    \end {document}